\documentclass[12pt,russian]{article}
\usepackage[english]{babel}
\usepackage[koi8-r]{inputenc} 
\usepackage{epsfig,latexsym,amsfonts,amsmath,amssymb}
\newcommand{\R}{\mathbb{R}}
\def\tc{\mathaccent"7017}

\def\XXint#1#2#3{{\setbox0=\hbox{$#1{#2#3}{\int}$}
	\vcenter{\hbox{$#2#3$}}\kern-.5\wd0}}

\begin{document}

\title{On symmetry and asymmetry\\
in a problem of shape optimization}
\author { {\it Alexander I. Nazarov}\footnote {Supported by RFBR grant 11-01-00825},\\ 
St.Peterburg Dept. of Steklov Institute,\\
St.-Petersburg State University, \\
{\small e-mail:\ al.il.nazarov@gmail.com}} 
\date{}\maketitle

\section{Introduction and main result}

Let $\Omega$ be an open subset of finite measure in $\R^n$. For
$1< p, q,r < \infty$ we consider the following extremal problem:
\begin{equation}\label{problem}
\lambda(p,q,r;\Omega)=\inf\Big\{\frac{\| \nabla u\|_{L_p(\Omega)}}
{\|u\|_{L_q(\Omega)}}: u\not\equiv0,\ u\in\tc{W}^1_p(\Omega),\ 
\int\limits_{\Omega}|u|^{r-2}u=0\Big\}.  
 \end{equation}
For $p=q=r=2$ the problem (\ref{problem}) is called
{\it twisted} Dirichlet eigenvalue problem, see \cite{BB}.
We are interested in the
so-called shape optimization problem: given $p,q,r$, among $\Omega$ with fixed volume, 
to find the set which minimizes $\lambda(p,q,r;\Omega)$. In other words our problem is 
a type of isoperimetric inequality. 

Even in one-dimensional case the problem (\ref{problem}) turned out to be rather complicated.
In this case the problem is to define the sets of positivity and negativity of the minimizer
to the problem
\begin{equation}\label{onedim}
\min\Big\{\frac{\Vert u^{\prime}\Vert _{L_p\left[-1,1\right]}}
{\Vert u \Vert_{L_q\left[-1,1\right]}}: u\left(\pm 1\right) =0,\ 
\int\limits_{-1}^1|u|^{r-2}u=0\Big\}.
\end{equation}
It was mentioned in \cite{DGS} that for $r=q$ the minimizer is always odd and thus the 
positivity and negativity sets are equal. On the other hand, authors of \cite{DGS} 
discovered that in the Poincar\'{e} case $r=2$ the situation is much more interesting: 
for $q\le 2p$ the minimizer is also odd while for $q>>1$ the symmetry 
breaking arises, and thus optimal partition of the interval $(-1,1)$ is not equal.

After efforts of many mathematicians (\cite{E}, \cite{BeKa}, \cite{Ka2}, \cite{A}) the final 
results in the problem (\ref{onedim}) were established in \cite{BKN}, \cite{N} for
$r=2$ and in \cite{DC} and recent paper \cite{GN} for general $r$. Namely, for
$q\le(2r-1)p$ the minimizer of (\ref{onedim}) is odd while $(2r-1)p$ is the bifurcation
point, and for $q>(2r-1)p$ the minimizer is asymmetric.

Note that the problem (\ref{onedim}) and equivalent statements in the Poincar\'{e} case arise 
also when estimating the eigenvalues in the Lagrange problem of the column stability 
(see \cite{EK} and references therein) and in the problem of local asymptotic optimality 
of non-parametric statistical tests (see \cite[Ch.6]{Ni}, \cite{LNN}).

\medskip
In multidimensional case the problem (\ref{problem}) and similar problems are considerably 
less investigated. We mention the papers \cite{GW} and \cite{N1} where symmetry and asymmetry 
of the minimizers were proved, respectively, in a ball and in a square, for $p=r=2$ and various 
values of $q$.

The next important point is the paper \cite{FH}, where the authors, dealing with the linear 
case $p=q=r=2$, by a very tight analysis solved the shape optimization problem and proved 
that the only optimal shape is given by a pair of equal balls.

In a quite recent paper \cite{CHP}, authors tried to generalize this results to the case
$r=q$ which is the simplest one in one-dimensional problem. They claimed that for any 
$p\in(1,\infty)$ and any admissible $q$ the optimal shape is also given by a pair of equal 
balls. In the present paper we show that it is not the case, and that, in contrast with
the problem (\ref{onedim}), the symmetry breaking phenomenon arises here.

\medskip
For any $p\in [1,\infty]$ we define the critical Sobolev exponent $p^*$ by relation
$$\frac 1{p^*}=\max\Big\{\frac 1p-\frac 1n;0\Big\}.
$$

Our main result is as follows.

\medskip
{\bf Theorem 1}. {\sl Let $1<p<\infty$ and $1<q<p^*$, $0<r-1<p^*$. Then

\begin{enumerate}
 \item The optimal shape in (\ref{problem}) is given only by a pair of disjoint balls.

\item For $q=r-1$ the optimal shape is given only by a pair of {\bf equal} balls.

\item For $q>\widehat q\equiv\big(\frac {r-1}n+1\big)^2p-\frac {(r-1)^2}n$ 
the optimal shape is given by a pair of {\bf non-equal} balls.

\end{enumerate}
}

{\bf Corollary 1}. For any $1<p<\infty$ and $1<r-1<p^*$, for sufficiently large
$q<p^*$ the optimal shape in (\ref{problem}) is asymmetric.

\medskip
{\bf Corollary 2}. The relation $q=r-1$ in Theorem 1, part 2, is the best possible in a sense.
Namely, given $\delta>0$, for any $1<p<n$ the optimal shape in (\ref{problem}) is asymmetric
provided $q=r-1+\delta$ is sufficiently close to $p^*$.

\medskip
The last corollary evidently disproves \cite[Theorem 1]{CHP}.

\medskip
{\bf Remark}. For $n=1$ Theorem 1, part 3, gives $\widehat q=r^2p-(r-1)^2$ that is greater than
the right exponent $(2r-1)p$, see \cite{DC} and \cite{GN}. To determine the true value of
the bifurcation parameter $\widehat q$ is an interesting open problem.

\section{Proof of Theorem 1}

Part 1 of Theorem 1 is quite standard. By compactness of embedding $\tc{W}^1_p(\Omega)$
into $L_q(\Omega)$ and into $L_{r-1}(\Omega)$, the infimum in (\ref{problem}) is attained.
Let $U$ be a minimizer. Define $U_{\pm}=\max\{\pm U,0\}$ and
$\Omega_{\pm}={\rm supp}(U_{\pm})$. Applying the Schwarz symmetrization 
(see, e.g., \cite{Ka}) to $U_{\pm}$ separately we obtain 
\begin{equation}\label{symm}
\lambda(p,q,r;\Omega)\ge\lambda(p,q,r;B_+\cup B_-),
\end{equation}
where $B_{\pm}$ are disjoint balls equimeasurable with $\Omega_{\pm}$, respectively. Moreover,
by the Euler equations both symmetrized functions $U_{\pm}^*$ have unique critical 
points. By \cite{BrZ}, in this case the numerator in (\ref{problem}) strictly decreases
under symmetrization, and thus the inequality in (\ref{symm}) is strict provided 
$\Omega\ne B_+\cup B_-$.

\medskip
To prove part 2, we write the Euler equation for (\ref{problem})
\begin{equation}\label{euler}
-{\rm div}(|\nabla U|^{p-2}\nabla U)=\lambda |U|^{q-2}U+\mu |U|^{r-2}
\end{equation}
(here $\lambda$ and $\mu$ are Lagrange's multipliers) and note that for $q=r-1$ we obtain
the following equations for $U_{\pm}$:
$$-{\rm div}(|\nabla U_{\pm}|^{p-2}\nabla U_{\pm})=(\lambda\pm\mu) |U_{\pm}|^{q-2}U_{\pm}.
$$
We know from part 1 that $\Omega= B_+\cup B_-$ and functions $U_{\pm}$ are radial.
Therefore, $U_{\pm}$ are in fact solutions of ODEs. By homogeneity, $U_{\pm}$ can be obtained
by dilation and multiplying by constant from a unique radial positive (generalized) solution of 
the boundary value problem
\begin{equation}\label{standard}
-{\rm div}(|\nabla v|^{p-2}\nabla v)=v^{q-1}\quad\mbox{in}\ \ B_1,\qquad v|_{\partial B_1}=0.
\end{equation}

Let $R_{\pm}$ be the radii of $B_{\pm}$. Then 
\begin{equation}\label{dilation}
U_{\pm}(x)=c_{\pm}\cdot v(\frac x{R_{\pm}}),
\end{equation}
and (\ref{problem}) is reduced to the following minimization problem:
\begin{equation}\label{finite}
\min\bigg\{\frac {c_+^pR_+^{n-p}+c_-^pR_-^{n-p}}{(c_+^qR_+^n+c_-^qR_-^n)^{\frac pq}} \ :
\ c_+^{r-1}R_+^n=c_-^{r-1}R_-^n,\ \ R_+^n+R_-^n=C\equiv\frac {|\Omega|}{|B_1|}.
\bigg\}
\end{equation}
Denoting $x=c_+^{r-1}R_+^n$ and $y=\frac 2C R_+^n-1\in(-1,1)$, we rewrite (\ref{finite}) 
as minimization of
\begin{equation}\label{one}
F(y)\equiv \frac {(1+y)^{1-\frac pn-\frac p{r-1}}+(1-y)^{1-\frac pn-\frac p{r-1}}}
{\big[(1+y)^{1-\frac q{r-1}}+(1-y)^{1-\frac q{r-1}}\big]^{\frac pq}}.
\end{equation}
Since $q=r-1$, the denominator in (\ref{one}) is constant. Since $q<p^*$, the exponent 
in the numerator is negative. Thus, the function (\ref{one}) is strictly convex, by symmetry 
its minimum is attained only at $y=0$, and the statement follows.

\medskip
To deal with part 3, we note that the radial positive solution of (\ref{standard}) solves
the problem
$$\inf\Big\{\frac{\| \nabla u\|_{L_p(B_1)}}
{\|u\|_{L_q(B_1)}}: u\not\equiv0,\ u\in\tc{W}^1_p(B_1)\Big\}.
$$
Therefore, if $R_+=R_-$ then the minimizer of (\ref{problem}) in $B_+\cup B_-$ satifies
(\ref{dilation}).

Now we consider functions (\ref{dilation}) for general $R_{\pm}$ and show that even such simple 
variations can provide symmetry breaking.

We again arrive at (\ref{one}) and write the Taylor series for the function $F$ at zero:
$$F(y)\sim2^{1-\frac pq}\cdot\Big(1+\frac {\gamma}2y^2\Big),\qquad
\gamma=\Big(\frac pn+\frac p{r-1}\Big)\Big(\frac pn+\frac p{r-1}-1\Big)
-\frac p{r-1}\Big(\frac q{r-1}-1\Big).
$$
For $q>\widehat q$ we have $\gamma<0$, and the statement follows.\hfill$\square$

\medskip
{\bf Proof of Corollaries 1 and 2}. For $p\ge n$ the statement of Corollary 1 is trivial.
Further, let $p<n$. By elementary calculation we obtain that the exponent
$\widehat q$ is increasing function of $r-1\in[0,p^*]$ and $\widehat q=p^*$ just for 
$r-1=p^*$. This gives both statements.\hfill$\square$


\begin{thebibliography}{AAAA}


\bibitem
{A}
J.M. Abessolo, 
{\em On a set of extremal problems and properties of related class of nonlinear
differential equations} 
// Ph.D. thesis, Moscow, 2002, 103p.

\bibitem
{BB}
L. Barbosa, P. B\'{e}rard, 
{\em Eigenvalue and ``twisted'' eigenvalue problems, applications to cmc surfaces} 
// J. Math. Pures Appl. (9) {\bf 79} (2000), N5, 427-450.

\bibitem
{BeKa}
M. Belloni, B. Kawohl, 
{\em A symmetry problem related to Wirtinger's and Poincare's inequality} 
// J. Diff. Eq. {\bf 156} (1999), 211-218.


\bibitem
{BrZ}
J.E. Brothers, W.P. Ziemer, 
{\em Minimal rearrangements of Sobolev functions} 
// J. Reine Angew. Math. {\bf 384} (1988), 153-179.


\bibitem
{BKN} 
A.P. Buslaev, V.A. Kondrat'ev, A.I. Nazarov, 
{\em On a family of extremal problems and related properties of an integral} 
// Mat. Zametki, {\bf 64} (1998), N6, 830-838 (Russian); 
English transl.: Math. Notes, {\bf 64} (1998), N5-6, 719-725. 

\bibitem
{DC}
G. Croce, B. Dacorogna, 
{\em On a generalized Wirtinger inequality} 
// Discr. Contin. Dyn. Syst. {\bf 9} (2003), N5, 1329-1341. 

\bibitem
{CHP}
G. Croce, A. Henrot, G. Pisante, 
{\em An isoperimetric inequality for a nonlinear eigenvalue problem} 
// Ann. Inst. H. Poincar\'{e} -- Anal. Nonlin., {\bf 29} (2012), 21-34.

\bibitem
{DGS}
B. Dacorogna, W. Gangbo, N. Subia, 
{\em Sur une generalisation de l'inegalite de Wirtinger} 
// Ann. Inst. H. Poincar\'e. Analyse Non Lin\'eaire. {\bf 9} (1992), 29-50.

\bibitem
{E}
Yu.V. Egorov, 
{\em On a Kondratiev problem} 
// C. R. A. S., Paris, Ser. I, {\bf 324} (1997), 503-507.

\bibitem
{EK}
Yu.V. Egorov, V.A. Kondratjev, 
{\em On a Lagrange problem} 
// C. R. A. S.,  Paris, Ser. I, {\bf 317} (1993), 903-918.     

\bibitem
{FH}
P. Freitas, A. Henrot, 
{\em On the first twisted Dirichlet eigenvalue} 
// Commun. Anal. Geom., {\bf 12} (2004), N5, 1083-1103.

\bibitem
{GN} 
I.V. Gerasimov, A.I. Nazarov, 
{\em Best constant in a three-parameter Poincar\'e inequality} 
// Probl. Mat. Anal. {\bf 61} (2011), 69-86 (Russian); 
English transl.: J. Math. Sci. {\bf 179} (2011), N1, 80-99. 


\bibitem
{GW} 
P. Gir\~{a}o, T. Weth, 
{\em The shape of extremal functions for Poincar\'{e}-Sobolev-type inequalities in a ball} 
// J. Func. An., {\bf 237} (2006), 194-223.

\bibitem
{Ka}
B. Kawohl, 
{\em Rearrangements and Convexity of Level Sets in PDE}, 
Lecture Notes in Mathematics, {\bf 1150}, Springer-Verlag, Berlin, 1985.


\bibitem
{Ka2}
B. Kawohl, {\em Symmetry results for functions yelding best constants in Sobolev type 
inequalities} 
// Discr. Contin. Dyn. Syst., {\bf 6} (2000), N3, 683-690.

\bibitem
{LNN}
M.A. Lifshits, A.I. Nazarov, Ya.Yu. Nikitin, 
{\em Tail behavior of anisotropic norms for Gaussian random fields} 
// C. R. A. S., Paris, Ser. I, {\bf 336} (2003), 85-88. 

\bibitem
{N}
A.I. Nazarov, 
{\em On exact constant in the generalized Poincar\'e inequality} 
// Probl. Mat. Anal., {\bf 24} (2002), 155-180 (Russian); 
English transl.: J. Math. Sci. {\bf 112} (2002), N1, 4029-4047. 

\bibitem
{N1}
A.I. Nazarov, 
{\em On the ``one-dimensionality'' of the extremal in the Poincar\'e inequality in the square}
// ZNS POMI, {\bf 259} (1999), 167-181 (Russian); 
English transl.: J. Math. Sci., {\bf 109} (2002), N5, 1928-1939. 

\bibitem
{Ni}
Ya.Yu. Nikitin, 
{\em Asymptotic Efficiency of Nonparametric Tests}, 
Cambridge University Press, 1995.







\end{thebibliography}
  \end{document}